\documentclass[12pt,reqno]{amsart}

\usepackage{latexsym}
\usepackage{graphicx}

\newtheorem{theorem}{{\sc Theorem}}[section]

\newtheorem{cor}[theorem]{{\sc Corollary}}

\newtheorem{lemma}[theorem]{{\sc Lemma}}

\theoremstyle{remark}

\def\f{\mathfrak }
\def\b{\mathbb }

\def\rr{\b R}
\def\mm{\b M}
\def\nn{\b N}
\def\cc{\b C}

\def\zz{\b Z}

\def\erw{\b E}
\def\hh{\b H}
\def\ii{\b I}
\def\jj{\b J}

\def\phi{\varphi }

\def\calf{{\mathcal F}}

\def\calc{{\mathcal C}}

\def\calp{{\mathcal P}}

\def\cals{{\mathcal S}}

\def\calx{{\mathcal X}}

\def\frg{{\f g}}

\def\frm{{\f m}}
\def\frn{{\f n}}

\def\frp{{\f p}}
\def\frq{{\f q}}

\def\on{\operatorname}

\def\tra{^{\prime}}

\begin{document}

\title[Wigner theorems]
{Wigner theorems for random matrices with dependent entries: Ensembles associated
to symmetric spaces and sample covariance matrices}

\author{Katrin Hofmann-Credner}
\address{Ruhr-Universit\"at Bochum, Fakult\"at f\"ur Mathematik, NA 3/31, D-44780 Bochum, Germany}
\email{katrin.hofmann-credner@ruhr-uni-bochum.de}
\author{Michael Stolz}
\address{Ruhr-Universit\"at Bochum, Fakult\"at f\"ur Mathematik, NA 4/32, D-44780 Bochum, Germany}
\email{michael.stolz@ruhr-uni-bochum.de}
\thanks{Research supported by Deutsche Forschungsgemeinschaft via SFB/TR 12
}

\date{\today}

\begin{abstract}
It is a classical result of Wigner that for an hermitian matrix
with independent entries on and above the diagonal, the mean
empirical eigenvalue distribution converges weakly to the
semicircle law as matrix size tends to infinity. In this paper, we
prove analogs of Wigner's theorem for random matrices taken from
all infinitesimal versions of classical symmetric spaces. This is
a class of models which contains those studied by Wigner and
Dyson, along with seven others arising in condensed matter
physics. Like Wigner's, our results are universal in that they
only depend on certain assumptions about the moments of the matrix
entries, but not on the specifics of their distributions. What is
more, we allow for a certain amount of dependence among the matrix
entries, in the spirit of a recent generalization of Wigner's
theorem, due to Schenker and Schulz-Baldes. As a byproduct, we
obtain a universality result for sample covariance matrices with
dependent entries.
\end{abstract}

\maketitle

\section{Introduction}

Classical physics-inspired random matrix theory is chiefly
concerned with probability measures on what Freeman Dyson in 1962
called the ``threefold way'', namely, the spaces of hermitian,
real symmetric, and quaternion real matrices (or their respective
exponentiated, compact versions). The rationale behind this focus
is Dyson's proof that any hermitian matrix (thought of as a
truncated Hamiltonian of a quantum system) that commutes with a
group of unitary symmetries and ``time reversals'' breaks down to
these three constituents (\cite{Dyson}), which are, in structural
terms, the tangent spaces to the Riemannian Symmetric Spaces (RSS)
of type A, AI and AII.

During the last decade, theoretical condensed matter physicists
have pointed out that matrix descriptions of systems such as
mesoscopic normal-superconducting hybrid structures are outside
the scope of Dyson's theorem, and that the tangent spaces to all
ten infinite series of classical RSS may (and do indeed) arise.
The deeper reasons are explained in \cite[Section 6.4]{AS},
\cite{AZ}, and \cite{HHZ}. Some of this material is summarized in
\cite{azldp}. Concrete matrix realizations of this ``tenfold way''
of ``symmetry classes'' are given in Section \ref{notat} below.

The task of developing random matrix theories for the full
``tenfold way'', i.e., studying probability measures on all ten
series of matrix spaces, has been taken up in \cite{azldp}, where
the probability measures enjoy invariance properties that
guarantee an explicit analytic expression for the joint eigenvalue
density, and in \cite{CSt, Duenez}, where the focus is on the
compact versions of the classical RSS, endowed with their natural
invariant probability measure. In the present paper, we abandon
invariance properties and turn to analogs of Wigner's famous
result of 1958 (\cite{Wigner58}), stating that for a symmetric
matrix with independent entries on and above the diagonal, the
mean empirical eigenvalue distribution converges weakly to the
semicircle law as matrix size tends to infinity. This is a
universality result in the sense that it only depends on certain
assumptions about the moments of the matrix entries, but not on
the specifics of their distributions.

Actually, our starting point is not the classical version of
Wigner's result, but a recent generalization, due to Schenker and
Schulz-Baldes (\cite{SSB}), allowing for a certain amount of
dependence to hold among the matrix entries. Specifically, the
authors consider the following set-up: For each $n \in \nn$ write
$\ii_n := \{ 1, \ldots, n\}$ and suppose that $\ii_n^2 = \ii_n
\times \ii_n$ comes with an equivalence relation $\sim_n$. The
entries of the matrix $X_n = ( \frac{1}{\sqrt{n}}\ a_n(p, q))_{p,
q = 1, \ldots, n}$ are complex random variables, with $a_n(p_1,
q_1), \ldots, a_n(p_j, q_j)$ independent whenever $(p_1, q_1),
\ldots, (p_j, q_j)$ belong to $j$ distinct equivalence classes of
the relation $\sim_n$. Furthermore, it is required that $a_n(p, q)
= \overline{a_n(q, p)}$ for all $n, p, q$. In the case that all
equivalence classes of $\sim_n$ are of the form $\{ (p, q), (q, p)
\}$, one is back to hermitian matrices with independent entries on
and above the diagonal, i.e.\ to the situation of Wigner's
theorem. If some equivalence classes are larger, then there is
some leeway for violations of independence. In the framework of
Schenker and Schulz-Baldes, the following conditions on $\sim_n$
serve as a less restrictive substitute for independence:

\begin{itemize}

\item[(W1)] $\max\limits_p \# \{(q,p\tra,q\tra)\in \ii_n^3:\
(p,q)\sim_n(p\tra, q\tra)\} = o(n^2)$

\item[(W2)] $\max\limits_{p,q,p\tra} \#\{q\tra\in \ii_n:\
(p,q)\sim_n(p\tra, q\tra)\} \leq B$, where $B<\infty$ is a
constant

\item[(W3)]$\#\{(p,q,p\tra)\in \ii_n^3:\
(p,q)\sim_n(q,p\tra)\mbox{ and } p\not= p\tra\} = o(n^2).$

\end{itemize}

Apart from that, one requires that for all $n, p, q$, $a_n(p, q)$
is centered and
\begin{equation}
\label{ssbvar} \erw ( a_n(p, q) \overline{a_n(p, q)}) = 1.
\end{equation}
Furthermore, a uniform bound on the $k$-th moments is assumed:
\begin{equation}
 \label{mobo} \sup_n\ \max_{p,q = 1, \ldots, n} \erw(
|a_n(p, q)|^k) < \infty\ \mbox{for all}\ k \in \nn.
\end{equation}

For an hermitian matrix $M \in \cc^{n \times n}$ with eigenvalues $\lambda_1,
\ldots, \lambda_n$ write
\begin{equation}
\label{Ln} L_n(M) := \frac{1}{n} \sum_{j=1}^n \delta_{\lambda_j}
\end{equation}
for the empirical measure of the eigenvalues of $M$. Then the main
theorem of \cite{SSB} can be stated as follows:
\begin{theorem}
\label{ssbthm} Under conditions (W1), (W2), (W3), \eqref{ssbvar},
and \eqref{mobo}, the measure $\erw(L_n(X_n))$, i.e.\ the mean
empirical eigenvalue distribution, converges weakly to the
semicircle law with density
\begin{equation}
\label{sc} \frac{1}{2\pi}\ 1_{[-2, 2]}(x)\   \sqrt{4-x^2}.
\end{equation}
\end{theorem}

So one ends up with the same limit distribution as in the
independent case. In fact, conditions (W1), (W2), (W3) arose from
a close reading of Wigner's proof, in order to understand how much
independence is really needed to arrive at the semicircle law. So
the approach of Schenker and Schulz-Baldes is complementary to the
route taken in a recent preprint of Anderson and Zeitouni
(\cite{AZeit}), where a different dependence structure leads to a
limit which is not a semicircle.

To establish analogs of Theorem \ref{ssbthm} for all ten symmetry
classes, we proceed as follows: In Section \ref{notat} we give
precise definitions of the matrix spaces in question and introduce
some auxiliary notation for the combinatorics of moment
calculations. In Section \ref{wdbdg} we treat those classes for
which $\erw L_n$ converges to the semicircle distribution and for
which nothing more than a slight extension of Theorem \ref{ssbthm}
is needed. On the other hand, substantial work has to be done for
the so-called ``chiral'' classes, which, despite their roots in
physics, are related to sample covariance matrices and hence lead
to some relative of the Mar\v{c}enko-Pastur distribution as limit
for $\erw L_n$. The main step is to rework the combinatorics of
moment convergence to the Mar\v{c}enko-Pastur law in the spirit of
Schenker and Schulz-Baldes, yielding a universality result for
sample covariance matrices with dependent entries. This is the
content of Section \ref{sample}. In the final Section
\ref{chiral}, this result is applied to the chiral classes.

\section{Background and notation}
\label{notat} We begin by listing the ten ``symmetry classes'',
i.e. series of matrix spaces, from which our matrices are taken.
In structural terms, these spaces are of the form $i \frg$, where
$\frg$ is the Lie algebra of a compact classical group, or
$i\frp$, where $\frp$ is the $-1$ eigenspace of a Cartan
involution of $\frg$ (see \cite{azldp} for details). The $i$
factor is to make sure that the matrices are hermitian. The labels
A, AI, etc.\ in the list are those of Lie theory, but no Lie
theoretic properties will be needed in what follows. $X^*$ denotes
the conjugate transpose of a matrix $X$.

\begin{description}

\item[Class A] $$  \mm_n^{{\rm A}} = \{ X \in {\mathbb C}^{n
\times n}:\ X\ \text{\rm hermitian} \}$$

\item[Class AI] $$  \mm_n^{{\rm AI}} = \{ X \in {\mathbb R}^{n
\times n}:\ X\ \text{\rm symmetric} \}$$

\item[Class AII] $$  \mm_n^{{\rm AII}} = \left\{ \left(
\begin{array}{rr} X_1 & X_2 \\ -\overline{X_2} & \overline{X_1}
\end{array} \right):\quad  \begin{array}{l} X_i \in {\mathbb C}^{n \times n}, X_1\ \text{\rm
hermitian},\\ X_2\ \text{\rm skew symmetric} \end{array} \right\}$$

\item[Class AIII] $$  \mm_n^{{\rm AIII}} = \left\{ \left(
\begin{array}{ll} 0 & X \\ X^* & 0 \end{array} \right):\ X \in {\mathbb C}^{s \times t} \right\}$$

\item[Class B/D] $$ \mm_n^{{\rm B/D}} = \{ X \in (i {\mathbb
R})^{n \times n}:\ X\ \text{\rm skew symmetric} \}$$

\item[Class BDI] $$  \mm_n^{{\rm BDI}} = \left\{ \left(
\begin{array}{ll} 0 & X \\ X^* & 0 \end{array} \right):\ X \in (i {\mathbb R})^{s \times t} \right\}$$

\item[Class DIII] $$ \mm_n^{{\rm DIII}} = \left\{ \left(
\begin{array}{rr} X_1 & X_2 \\ X_2 & -X_1
\end{array} \right):\ X_i \in (i {\mathbb R})^{n \times n}\ \text{\rm skew symmetric} \right\}$$

\item[Class C] $$ \mm_n^{{\rm C}} = \left\{ \left(
\begin{array}{rr} X_1 & X_2 \\ \overline{X_2} & - \overline{X_1}
\end{array} \right):\quad \begin{array}{l} X_i \in \cc^{n \times n},\ X_1\ {\rm hermitian},\\ X_2\ {\rm symmetric}\end{array}  \right\}$$

\item[Class CI] $$  \mm_n^{{\rm CI}} = \left\{ \left(
\begin{array}{rr} X_1 & X_2 \\ X_2 & -X_1
\end{array} \right):\ X_i \in {\mathbb R}^{n \times n}\ \text{\rm symmetric} \right\}$$

\item[Class CII] $$  \mm_n^{{\rm CII}} = \left\{ \left(
\begin{array}{ll} 0 & X \\ X^* & 0 \end{array} \right):\ X \in {\mathbb H}^{s \times t} \right\},$$
where the space $\hh^{s \times t}$ of quaternionic matrices is
embedded into $\cc^{2s \times 2t}$ as
$$ \hh^{s \times t} = \left\{ \left( \begin{array}{rr} U & V \\ - \overline{V} & \overline{U}\end{array}\right):\ U, V \in \cc^{s\times t}\right\}.$$

\end{description}

Generically, we write $\calc$ for any label A, AI,..., CII. We
have $\mm_n^{\calc} \subset \cc^{\delta n \times \delta n}$ with
$\delta = \delta_{\calc} = 2$ if $\calc =$ AII, DIII, C, CI, CII
and $\delta = 1$ otherwise. Classes AIII, BDI and CII, the
``chiral'' classes in physics terminology (since they are related
to Dirac operators, see \cite{Forr}, \cite{HHZ}), are special in
that the shape of the subblocks depends on an extra parameter $s$,
and it is clear that one will have to control the relative growth
of $s = s(n)$ and $t(n) = n - s(n)$ as $n \to \infty$. In fact,
their large $n$ behaviour is quite different from that of the
other classes. Consequently, in what follows, we will treat these
classes separately. A, AI, and AII are the classical Wigner-Dyson
classes (underlying GUE, GOE, and GSE, resp.). The remaining
classes, arising via Bogolioubov-de Gennes mean field
approximation (see \cite[Sec. 6.4]{AS}), will be termed
Bogolioubov-de Gennes (BdG) or superconductor classes.

 Before getting down to business, let us review some
combinatorial notation and facts that will be useful later on.
Write $\calp(n)$ for the set of all partitions of $\ii_n = \{ 1,
2, \ldots, n\}$. If $\frp \in \calp(n)$ has $r$ blocks, write
$|\frp| = r$. Define $\calp^{(i)}(n) := \{ \frp \in \calp(n):\
|\frp| = i\}$. If each of the blocks of $\frp$ consists of exactly
two elements, we say that $\frp$ is a pair partition and write
$\frp \in \calp_2(n)$. For $\frp \in \calp(n)$, write
$\sim_{\frp}\ \subseteq\ \ii_n \times \ii_n$ for the corresponding
equivalence relation. Given $\frp, \frq \in \calp(n)$, define
$\sim_{\frp} \vee \sim_{\frq}\ \subseteq\ \ii_n \times \ii_n$ as
the smallest equivalence relation which contains $\sim_{\frp}$ and
$\sim_{\frq}$. The partition corresponding to $\sim_{\frp} \vee
\sim_{\frq}$ is denoted by $\frp \vee \frq$.

$\frp \in \calp(n)$ is called crossing, if there exist $p_1 < q_1
< p_2 < q_2$ in $\ii_n$ such that $p_1 \sim_{\frp} p_2
\not\sim_{\frp} q_1 \sim_{\frp} q_2$. Otherwise, it is called
noncrossing. Write ${\rm NC}(n),\ {\rm NC}^{(i)}(n)$ and ${\rm
NC}_2(n)$ for the set of noncrossing partitions, noncrossing
partitions with $i$ blocks, and noncrossing pair partitions of
$\ii_n$, respectively. For sets $\Omega, \Omega\tra$ write $
\calf(\Omega, \Omega\tra) := \{ \phi: \Omega \to \Omega\tra\}$ and
$\calf(k, n) := \calf(\ii_k, \ii_n).$ For $\frp \in \calp(k)$
write $\calf(\frp, \Omega)$ for the set of all $\phi \in
\calf(\ii_k, \Omega)$ which are constant on the blocks of $\frp$.
Finally, let us give names to two important special pair
partitions:
\begin{eqnarray*}
\frm &:=& \{ \{ 1, 2\}, \{3, 4\}, \ldots, \{ 2k-1, 2k\}\}\in
\calp_2(2k)\quad
\text{and}\\
\frn &:=& \{ \{ 2, 3\}, \{ 4, 5\}, \ldots, \{ 2k, 1\} \}\in
\calp_2(2k).
\end{eqnarray*}

It is well-known that $\# \on{NC}_2(2k)$ equals the $k$th Catalan
number $C_k$, which in turn equals the $2k$th moment of the
semicircular distribution with density given in \eqref{sc} (whose
odd moments vanish). On the other hand, setting for $\kappa > 0$
\begin{equation}
\label{mommp} m_k := \sum_{i=1}^k \# (\on{NC}^{(i)}(k))\
\kappa^{i},
\end{equation}
$(m_k)_{k \in \nn}$ is the sequence of moments of the
Mar\v{c}enko-Pastur distribution with density
$$
\max\{0, (1 - \kappa)\} \delta_0 + \frac{ \sqrt{ 4 \kappa -
(x-1-\kappa)^2}}{2\pi x}\ 1_{[(1- \sqrt{\kappa})^2,\ (1+
\sqrt{\kappa})^2]}(x).
$$
A reference for these facts is \cite{HT} or \cite{HP}.

\section{Wigner-Dyson and Bogolioubov-de Gennes classes}
\label{wdbdg}

These are the easy cases, because they basically reduce to Theorem
\ref{ssbthm}. In fact, one may interpret the symmetries of a
matrix from $\mm_n^{\calc}$ as an equivalence relation on pairs of
indices, recalling that in the set-up of Theorem \ref{ssbthm}, if
index pairs $(p, q), (p\tra, q\tra)$ are equivalent, then the
corresponding random variables $a_n(p, q),\ a_n(p\tra, q\tra)$ may
be identical. Of course, for $\calc =$ AII, DIII, C, CI, CII, the
symmetries of $\mm_n^{\calc}$ must be realized as equivalence
relations on $\ii_{2n} \times \ii_{2n}$, so one obtains the
desired theorem on random elements of $\mm_n^{\calc}$ by passage
to a subsequence in Theorem \ref{ssbthm}. A more serious caveat is
the following: Some of the blocks which make up matrices from
$\mm_n^{\calc}$ are skew symmetric, so their diagonal elements
vanish, contradicting condition \eqref{ssbvar}. While we will see
that this problem can be circumvented in the cases at hand, the
full blocks of zeroes in the chiral cases make it impossible to
apply Theorem \ref{ssbthm} for them as well.

To make the set-up for this section precise, let $\calc$ be a
Wigner-Dyson or BdG class, write $\delta = \delta^{\calc}$ as in
Section \ref{notat}, $\jj_{\delta n}^{\calc} := \{ (p, q) \in
\ii_{\delta n}^2:\ \operatorname{pr}_{p, q}(\mm_n^{\calc}) \neq
0\}$, where $\operatorname{pr}_{p, q}$ projects each element of
$\mm_n^{\calc}$ onto its $(p, q)$-entry. Consider an equivalence
relation $\sim_{\delta n}$ on $\ii_{\delta n}^2$ and a random
matrix  $X_{\delta n} = ( \frac{1}{\sqrt{\delta n}}\ a_{\delta
n}(p, q))_{p, q = 1, \ldots, \delta n}$ such that the centered
complex random variables $a_{\delta n}(p_1, q_1), \ldots,
a_{\delta n}(p_j, q_j)$ are independent whenever $(p_1, q_1),
\ldots, (p_j, q_j)$ belong to $j$ distinct equivalence classes of
the relation $\sim_{\delta n}$. Assume that conditions (W1), (W2),
(W3) hold, with $\ii_{\delta n}$ in the place of $\ii_n$, and that
the moment condition \eqref{mobo} is satisfied. As to
\eqref{ssbvar}, it is required that it holds for all $(p, q) \in
\jj_{\delta n}^{\calc}$. All realizations of the matrix $X_{\delta
n}$ are supposed to be elements of $\mm_n^{\calc}$. It is
straightforward to verify that this assumption is compatible with
conditions (W1), (W2), (W3). So we may take the symmetries of the
matrices for granted, and have some leeway for extra dependence
between the matrix entries. Under these conditions, there holds
\begin{theorem}
\label{thmwdbdg} If $\calc$ is a Wigner-Dyson or Bogolioubov-de
Gennes class, $\erw(L_{\delta n}(X_{\delta n}))$ converges weakly to the
semicircle law.
\end{theorem}
It only remains to address the
complication that random elements of $\mm_{n}$ may have up
to $4n$ entries which are identically zero. To see that the effect
of this complication is asymptotically negligible, recall from the proof of Theorem \ref{ssbthm} in
\cite{SSB} that the $k$-th moment of $\erw L_n$ vanishes if $k$ is
odd and is asymptotically equivalent to
$$ \frac{1}{n^{l + 1}} \sum_{\frp \in \on{NC}_2(2l)}
\# \cals_n(\frp)$$ if $k= 2l$ is even. Here the set
$\cals_n(\frp)$ consists of all pairs $(\phi, \psi),\ \phi, \psi
\in \calf(k, n)$, with the following properties:
\begin{itemize}
\item[(i)] $\psi(j) = \phi(j+1)$ for all $j \in \ii_k$, where $k +
1$ is cyclically identified with $1$.

\item[(ii)]
\begin{eqnarray*}
( \phi(i), \psi(i)) = (\psi(j), \phi(j)) &\text{if}& i \sim_{\frp}
j\quad \text{and}\\
( \phi(i), \psi(i)) \not\sim_n (\phi(j), \psi(j)) &\text{if}& i
\not\sim_{\frp} j.
\end{eqnarray*}

\end{itemize}

For $n \in \nn$ fix $E_n \subset \ii_{n}^2$ with $\# E_n =
o(n^2)$. Actually, what we have in mind is that $E_n$ contains the
$O(n)$ diagonal places of skew blocks. For $\nu \in \ii_k$ set
$$ \cals_n^{(\nu)}(\frp) := \{ (\phi, \psi) \in \cals_n(\frp):\
(\phi(\nu), \psi(\nu)) \in E_n\}$$ and
$$ \cals_n(\frp)\tra = \bigcup_{\nu \in \ii_k}
\cals_n^{(\nu)}(\frp).$$ Then the following lemma makes it
possible to neglect the effect of the zero entries on the
diagonals of the blocks:
\begin{lemma}
\label{nullen} For $\frp \in \on{NC}_2(2l)$,\ $\#
\cals_n(\frp)\tra = o(n^{l+1}).$
\end{lemma}
\begin{proof}
Since  $ \# \cals_n(\frp)\tra \le \sum_{\nu = 1}^{2l} \#
\cals_n^{(\nu)}(\frp),$ it suffices to show that for all $\nu$ one
has $\# \cals_n^{(\nu)}(\frp) = o(n^{l+1}).$ To this end, we
construct an element $(\phi, \psi) \in \cals_n^{(\nu)}(\frp)$,
starting with $(\phi(\nu), \psi(\nu))$ and proceeding cyclically.
By cyclically permuting the index set, we may assume that $\nu =
1$. For $(\phi(1), \psi(1))$ we have $o(n^2)$ choices. $\phi(2)$
is then fixed by property (i). If $1 \sim_{\frp} 2$, then
$\psi(2)$ is fixed by (ii). Otherwise, we have at most $n$
choices. Similarly, for $j \ge 2$, once we have chosen
$\psi(j-1)$, $\phi(j)$ is fixed, and $\psi(j)$ is either fixed or
we have at most $n$ choices for it, depending on whether $i
\sim_{\frp} j$ for some $1 \le i < j$ or not. The latter case
occurs $l-1$ times. So $\# \cals_n^{(\nu)}(\frp) \le o(n^2) \times
n^{l-1} = o(n^{l+1}).$

\end{proof}

\section{Sample covariance matrices}
\label{sample} In this section we prove a limit theorem for the
mean empirical eigenvalue distribution of sample covariance
matrices with dependence. By these we understand matrices of the
form $X^*X$, where $X$ admits a certain amount of dependence among
its entries, $X^*$ is the conjugate transpose of $X$, and the
entries of $X$ are not necessarily Gaussian, but subject to
certain conditions on their moments. In our context, this is
preparatory work for the study of the chiral classes in Section
\ref{chiral}, but it is of interest for its own sake. In the case
that $X$ has independent entries, the result is well-known, and we
will take the combinatorial proof of Oravecz and Petz (\cite{OP})
as starting-point for an analysis in the spirit of Schenker and
Schulz-Baldes (\cite{SSB}).

For $n \in \nn$ let $s(n), t(n) \in \nn$ and suppose that there
exist $\kappa, \mu \in ]0, \infty[$ such that $\lim_n
\frac{s(n)}{n} = \kappa$ and $\lim_n \frac{t(n)}{n} = \mu$. The
classical case is $t(n) = n$, but we will need this more general
set-up in Section \ref{chiral}. Consider an equivalence relation
$\sim_n$ on $\ii_{s(n)} \times \ii_{t(n)}$ and a random matrix
$X_n = ( \frac{1}{\sqrt{n}}\ a_n(p, q))_{p = 1, \ldots, s(n), q =
1, \ldots, t(n)}$ such that the complex random variables $a_n(p_1,
q_1), \ldots, a_n(p_j, q_j)$ are independent whenever $(p_1, q_1),
\ldots, (p_j, q_j)$ belong to $j$ distinct equivalence classes of
the relation $\sim_n$. We impose the following conditions on
$\sim_n$:

\begin{itemize}

\item[(MP1)] $\max\limits_p \# \{(q,p\tra,q\tra)\in \ii_{t(n)}
\times \ii_{s(n)} \times \ii_{t(n)}:\ (p,q)\sim_n(p\tra, q\tra)\}
= o(n^2).$

\item[(MP2)] $\max\limits_{p,q,p\tra} \#\{q\tra\in \ii_{t(n)}:\
(p,q)\sim_n(p\tra, q\tra)\} \vee \max\limits_{p,q, q\tra}
\#\{p\tra\in \ii_{s(n)}:\ (p,q)\sim_n(p\tra, q\tra)\} \le B$,
where $B$ is a finite constant.

\item[(MP3)] $\#\{(p,q, q\tra)\in \ii_{s(n)} \times \ii_{t(n)}^2
:\ (p,q)\sim_n(p, q\tra)\mbox{ and } q\neq q\tra\} = o(n^2)$ and\\
$\#\{(p, p\tra, q)\in \ii_{s(n)}^2 \times \ii_{t(n)} :\
(p,q)\sim_n(p\tra,q)\mbox{ and } p \neq p\tra\} = o(n^2)$.
\end{itemize}
We assume that \eqref{ssbvar} and \eqref{mobo} hold. Under these
assumptions we will prove the following theorem.

\begin{theorem}
\label{thmwish} As $n \to \infty$, $\erw(L_n(X_n^* X_n))$
converges weakly to a probability measure with $k$-th moment equal
to
\begin{equation}
\label{genmom}
 \sum_{i=1}^k \# (\on{NC}^{(i)}(k))\ \kappa^{i} \mu^{k - i + 1}.
\end{equation}
If $\mu =1$, this limit is the Mar\v{c}enko-Pastur distribution.
\end{theorem}

We are going prove Theorem \ref{thmwish} via the method of
moments. So we fix $k \in \nn$ and show that
\begin{equation}
\label{lnmom} \int x^k \erw(L_n(X_n^* X_n))(dx) = \erw \int x^k
L_n(X_n^* X_n)(dx) = \frac{1}{n}\ \erw \on{Tr}((X_n^* X_n)^k)
\end{equation}
converges to \eqref{genmom} as $n \to \infty$. To write the trace
in \eqref{lnmom} in an explicit way, we use the notation
introduced in Section \ref{notat}. In the course of the technical
proofs it will be convenient to identify the index set $\ii_{2k}$
with the cyclic group $\zz / 2k\zz$, i.e., to identify $2k + 1$
with $1$ and so on.

The starting point for all that follows is the observation that

\begin{equation}
\label{qv1} \frac{1}{n} \on{Tr}((X_n^* X_n)^k) = \frac{1}{n^{k+1}}
\sum_{\phi \in \calf(\frm, s(n))} \sum_{\psi \in \calf(\frn,
t(n))}  \prod_{\nu = 1}^k\ \overline{a_n( \phi(2\nu - 1),
\psi(2\nu - 1))}\ a_n(\phi(2\nu), \psi(2\nu)),
\end{equation}
with $\frm, \frn$ as in Section \ref{notat}.
Let $\frp \in \calp(2k)$. We say that $(\phi, \psi) \in
\calf(\frm, s(n)) \times \calf(\frn, t(n))$ is associated to
$\frp$, and write $(\phi, \psi) \in S_n(\frp)$, if for all $i, j =
1, \ldots, 2k$ there holds
\begin{equation}
\label{assoc} i \sim_{\frp} j\ \Longleftrightarrow\ (\phi(i),
\psi(i)) \sim_n (\phi(j), \psi(j)).
\end{equation}
Writing
\begin{equation}
\label{qv4}
 \Sigma_n(\frp) := \sum_{(\phi, \psi) \in S_n(\frp)}
\prod_{\nu = 1}^k\ \overline{a_n( \phi(2\nu - 1), \psi(2\nu -
1))}\ a_n(\phi(2\nu), \psi(2\nu)),\end{equation}
 we obtain
\begin{equation}
\label{qv2} \frac{1}{n} \on{Tr}((X_n^* X_n)^k) = \frac{1}{n^{k+1}}\
\sum_{\frp \in \calp(2k)}\ \Sigma_n(\frp).
\end{equation}

\begin{lemma}
\label{basic} For $\frp \in \calp(2k)$ one has
$\# S_n(\frp) \le O(n^{|\frp| + 1}).$
\end{lemma}
\begin{proof}
We construct an element $(\phi, \psi) \in S_n(\frp)$, proceeding from $1$ to $2k$ and giving rather coarse upper bounds on the number of choices in each step.
We have $s(n)$ choices
for $\phi(1)$ and $t(n)$ choices for $\psi(1)$. For $\phi$ to be constant on the blocks of $\frm$,
we must have $\phi(2) = \phi(1)$. If $1 \sim_{\frp} 2$, then by (MP2) we have at most $B$
choices for $\psi(2)$. Otherwise, we have at most $t(n)$ choices.
Note that $\psi(2) = \psi(3)$, since
$\psi$ is supposed to be constant on the blocks of $\frn$. In the general case, for $\nu = 2, 3, \ldots, 2k-1$, one of $\phi(\nu), \psi(\nu)$
is fixed, and for the other we have $\le B$ choices if $\nu \sim_{\frp} \nu\tra$ for some $\nu\tra \in \ii_{\nu - 1}$ or
at most $s(n) \vee t(n)$ choices otherwise. This latter case occurs precisely $|\frp| - 1$ times.
$(\phi(2k), \psi(2k))$ is fixed by the requirement that $\phi(2k-1) = \phi(2k)$ and $\psi(2k) = \psi(1)$. In total, we have at most
$s(n)\ t(n)\ (s(n)\vee t(n))^{|\frp| - 1}\ B^{2k - (|\frp| + 1)} = O(n^{|\frp| + 1})$ choices for $(\phi, \psi)$.
\end{proof}

\begin{lemma}
\label{ck} For $\frp \in \calp(2k)$, $|\erw \Sigma_n(\frp)| \le \# S_n(\frp) c_k$, with $c_k$ independent of $n$.
\end{lemma}
\begin{proof}
By H\"older's inequality and \eqref{mobo} one has
\begin{eqnarray*} &&\left| \erw \prod_{\nu = 1}^k\ \overline{a_n(
\phi(2\nu - 1), \psi(2\nu - 1))}\ a_n(\phi(2\nu),
\psi(2\nu))\right| \\   &\le& \prod_{\nu =
1}^{2k}\ \  \erw \left( \left| a_n( \phi(\nu), \psi(\nu
))\right|^k\right)^{\frac{1}{k}} \le c_k < \infty,
\end{eqnarray*}
with $c_k$ independent of $n, \phi$ and $\psi$.
\end{proof}

\begin{cor}
\label{pair}
$\erw \Sigma_n(\frp) = o(n^{k+1})$ unless $\frp \in \calp_2(2k)$.
\end{cor}
\begin{proof}
If
$|\frp| \le k - 1$, then by Lemmata \ref{basic} and \ref{ck},
$ | \erw \Sigma_n(\frp)| \le \# S_n(\frp)\ c_k \le O(n^k) c_k =
o(n^{k+1}).$
If $\frp$
contains a block which consists of precisely one element $\nu_0$,
say, then we have $\erw \Sigma_n(\frp)= 0$, because for any
$(\phi, \psi) \in S_n(\frp)$ the random variable $a_n(\phi(\nu_0),
\psi(\nu_0))$ is centered and by construction independent of $\{
a_n(\phi(\nu), \psi(\nu)):\ \nu \in \ii_{2k} \setminus \{
\nu_0\}\}.$ So $\erw \Sigma_n(\frp)$ vanishes if $|\frp| \ge k+1$, or if $|\frp| = k$, but $\frp \not\in \calp_2(2k).$
\end{proof}

The following lemma is a straightforward adaptation of a key step
of \cite{SSB} to the present context.

\begin{lemma}
\label{cr}
If $\frp \in \calp_2(2k)$ is crossing, then $\# S_n(\frp) =
o(n^{k+1})$.
\end{lemma}
\begin{proof}
Suppose that $\frp$ contains a nearest neighbour pair, i.e.\ a
block of the form $\{\nu, \nu+1\}$.  Assume that $\nu$ is odd.

If
$\psi(\nu) = \psi(\nu + 1)$, then $\psi(\nu - 1) = \psi(\nu + 2)$.
Writing $\jj = \ii_{2k} \setminus \{ \nu, \nu + 1\}$, we see that
$(\phi|_{\jj}, \psi|_{\jj}) \in S_n(\tilde{\frp})$, where
$\tilde{\frp}$ is the partition of $\jj$ whose blocks are those of
$\frp$ except for $\{ \nu, \nu + 1\}$. Given $(\tilde{\phi},
\tilde{\psi}) \in S_n(\tilde{\frp})$, there are $s(n)$ choices for
$\phi(\nu) = \phi(\nu + 1)$, hence $s(n)$ possible extensions to
$(\phi, \psi) \in S_n(\frp)$, since $\psi(\nu), \psi(\nu + 1)$ are
determined by $\tilde{\psi}$.

If $\psi(\nu) \neq \psi(\nu + 1)$,
then by (MP3) there are only $o(n^2)$ choices for the triplet
$(\phi(\nu) = \phi(\nu + 1), \psi(\nu), \psi(\nu + 1))$. As in the
proof of Lemma \ref{basic}, we see that there are at most
$O(n^{k-1})$ choices for the remaining values of $\phi, \psi$,
since $\tilde{\frp}$ consists of $k-1$ pairs.

Since we may argue
analogously for $\nu$ even, we have shown that
\begin{equation}
\label{qv10} \# S_n(\frp) \le O(n)\ \#S_n(\tilde{\frp}) +
o(n^{k+1}).
\end{equation}
Since $\frp$ was assumed to be crossing, iterating this argument
yields
\begin{equation}
\label{qv11}
 \# S_n(\frp) \le O(n^r)\ \#S_n(\frp\tra) + o(n^{k+1}),
\end{equation}
where $\frp\tra$ contains no nearest neighbour pair. Upon
relabelling, we may regard $\frp\tra$ as an element of
$\calp_2(2(k-r))$, where $k - r \ge 2$. Note that we may have
$\frp = \frp\tra$, whence it is possible that $r = 0$.

Let $\lambda$ be the minimal $l > 0$ such that there exists $\nu
\in \ii_{2(k-r)}$ with $\{ \nu, \nu+l\} \in \frp\tra$ (where
addition takes place in $\zz / 2(k-r) \zz$). Observe that $\lambda
\ge 2$. Now, if $\nu$ has the property that $\{ \nu, \nu + \lambda\} \in \frp\tra$, then all elements of $\{ \nu + 1,
\ldots, \nu + \lambda - 1\}$ are paired with partners outside $\{
\nu, \nu + 1, \ldots, \nu + \lambda\}$. We find an upper bound for
$\# S_n(\frp\tra)$ as follows. By (MP1), we have $s(n)$ choices
for $\phi(\nu)$ and $o(n^2)$ choices for the triplet $(\psi(\nu),
\phi(\nu + \lambda), \psi(\nu + \lambda))$. By construction, going
from left to right through $\nu + l\ (l = 1, 2, \ldots, \lambda -
2)$, either $\phi(\nu + l)$ is fixed and there are at most $t(n)$
choices for $\psi(\nu + l)$, or $\psi(\nu + l)$ is fixed and there
are at most $s(n)$ choices for $\phi(\nu + l)$.
According to whether $\nu$ is even or odd, we must have
either $\phi(\nu + \lambda - 1) = \phi(\nu + \lambda)$ and
$\psi(\nu + \lambda - 1) = \psi(\nu + \lambda - 2)$ or $\psi(\nu +
\lambda - 1) = \psi(\nu + \lambda)$ and $\phi(\nu + \lambda - 1) =
\phi(\nu + \lambda - 2)$. So we have at most $O(n^{\lambda - 1})
o(n^2)$ choices for the restrictions of $\phi$ and $\psi$ to $\{
\nu, \nu + 1, \ldots, \nu + \lambda\}$. Going cyclically through
the complement of this set, starting with $\nu + \lambda + 1$, in
each step one of $\phi(\nu + \lambda + l), \psi(\nu + \lambda +
l)$ is fixed, and there are at most $B$ resp.\ $O(n)$ choices for
the other, depending on whether $\nu + \lambda + l$ is paired with
one of the previously considered points or not. This latter case
occurs exactly $k - r - \lambda$ times.

Putting all this together with \eqref{qv11}, we arrive at
$$ \# S_n(\frp) \le O(n^r)\ o(n^{\lambda + 1})\ O(n^{k - r - \lambda}) +
o(n^{k+1}) = o(n^{k+1}).$$

\end{proof}

\begin{lemma}
\label{nc} For all $\frp \in \on{NC}_2(2k)$ there holds $l_{\frp}
:= |\frp \vee \frm| + |\frp \vee \frn| = k+1$, with $\vee$ as in
Section \ref{notat}.
\end{lemma}
\begin{proof}

For $k=1$, the only (noncrossing pair) partition $\frp \in
\calp(2)$ is $\{ \{ 1, 2\}\}$, which satisfies  $l_{\frp} = 2$.
Suppose that $k \ge 2$ and that the claim is true for $k-1$. Note
that any $\frp \in \on{NC}_2(2k)$ contains a block of the form $\{ \nu, \nu +
1\}$. Now consider the partition $\tilde{\frp}$ of $\ii_{2k}
\setminus \{ \nu, \nu + 1\}$, whose blocks are those of $\frp$
except for $\{ \nu, \nu + 1\}$. Define $\tilde{\frm}, \tilde{\frn}
\in \ii_{2k} \setminus \{ \nu, \nu + 1\}$ in the obvious way. We claim that
\begin{equation}
\label{lp2} l_{\tilde{\frp}} = l_{\frp} - 1.
\end{equation}
To see this for $\nu$ even, note that in this case $\{ \nu, \nu +1\}$ is
also a block of $\frn$. This means that it is a block of $\frp
\vee \frn$, hence $|\tilde{\frp} \vee \tilde{\frn}| = |\frp \vee
\frn| -1$. On the other hand, $\frp \vee \frm$ has a block which
contains $\nu - 1, \nu, \nu + 1, \nu + 2$. This block does not split on eliminating
$\{ \nu, \nu + 1\}$, since $\{ \nu - 1, \nu + 2\}$ is a block of $\tilde{\frm}$.
So we arrive at $|
\tilde{\frp} \vee \tilde{\frm}| = |\frp \vee \frm|$. For $\nu$
odd, the argument is analogous, yielding \eqref{lp2}.
 Hence
by induction, $l_{\frp} = l_{\tilde{\frp}} + 1 = ((k-1) + 1) + 1 =
k+1$.
\end{proof}

The following is evident:
\begin{lemma}
\label{evenodd} For $\frp \in \on{NC}_2(2k)$, each of the blocks
of $\frp$ consists of exactly one odd and exactly one even number.
\end{lemma}

\begin{lemma}
\label{count} $  \# \{ \frp \in \on{NC}_2(2k):\ |\frp \vee \frm|
= j\ \} = \# \on{NC}^{(j)}(k).$
\end{lemma}
\begin{proof}
We may identify $\calp(k)$ with $ \{ \frp \in \calp(2k):\
\text{any block of}\ \frp\ \text{is the union of blocks of}\ \frm
\}$. So $\frp \mapsto \frp \vee \frm$ maps $\calp(2k)$ onto
$\calp(k)$. It is easy to see that if $\frp \vee \frm$ is
crossing, so is $\frp$. In fact $\on{NC}_2(2k)$ is mapped
bijectively onto $\on{NC}(k)$. To see this, it suffices to show
that $\frp \mapsto \frp \vee \frm$ is injective, as it is known
that $\# \on{NC}(k) = \# \on{NC}_2(2k)$ (see \cite[Remark
9.5]{NSp}). A block of $\frp \vee \frm$ is of the form $b_J =
\bigcup_{\nu \in J} \{ 2\nu - 1, 2\nu\}$ for some $J \subseteq
\ii_k$. We have to show that there exists precisely one
$\tilde{\frp} \in \on{NC}_2(b_J)$ such that $\tilde{\frp} \vee \{
\{ 2\nu - 1, 2\nu\}:\ \nu \in J\} = \{ b_J\}$. Since our aim is to
show that a pairing of the elements of $b_J$ with certain
properties is uniquely determined, the embedding of $b_J$ into
$\ii_{2k}$ is irrelevant, and we may assume that $b_J = \ii_{2r}$
for some $r \le k$. Let us start by finding a partner for $1$. By
Lemma \ref{evenodd}, the partner must be even, $1
\sim_{\tilde{\frp}} 2s$, say. Assume that $s < r$. Since we wish
to construct a noncrossing $\tilde{\frp}$, no $1 < \nu < 2s$ can
be paired with any $\nu\tra > 2s$. On the other hand, $2s
\not\sim_{\frm} 2s+1$, so $2s$ and $2s+1$ lie in distinct blocks
of $\tilde{\frp} \vee \frm$, contradicting the requirement that
$\tilde{\frp} \vee \frm = \{ \ii_{2r}\}$. Consequently, we must
have $1 \sim_{\tilde{\frp}} 2r$.

The partner of $2$ must be odd. We claim that necessarily $2
\sim_{\tilde{\frp}} 3$. Otherwise the partner of $2$ is $2s - 1,\
s \ge 3$. Then for $\tilde{\frp}$ to be noncrossing, it is
necessary that $\{ 3, 4, \ldots, 2s-2)\}$ be a union of blocks of
$\tilde{\frp}$. But $ 2 \not\sim_{\frm} 3$ and $2s-2
\not\sim_{\frm} 2s-1$. So $\{ 3, 4, \ldots, 2s-2)\}$ splits into at least
two distinct blocks of $\tilde{\frp} \vee \frm$, violating our
conditions on $\tilde{\frp}$. Deleting $\{ 2, 3\}$ and
relabelling, we inductively see that $\tilde{\frp} = \{ \{ 1,
2r\}, \{2, 3\}, \{ 4, 5\}, \ldots, \{ 2r - 2, 2r - 1\}\}$, hence
is uniquely determined.
\end{proof}

Define
$$ S_n^{\vee}(\frp) := \{ (\phi, \psi):\ \phi \in
\calf(\frp \vee \frm, s(n)),\ \psi \in \calf(\frp \vee \frn,
t(n))\}$$ and observe that $S_n^{\vee}(\frp) \subset S_n(\frp)$.
Write
$$ S_n^{\wedge}(\frp) : = S_n(\frp) \setminus S_n^{\vee}(\frp)$$
and
$$ \Sigma_n^{\vee}(\frp) := \sum_{(\phi, \psi) \in S_n^{\vee}(\frp)} \prod_{\nu = 1}^k\ \overline{a_n( \phi(2\nu - 1),
\psi(2\nu - 1))}\ a_n(\phi(2\nu), \psi(2\nu)).$$

\begin{lemma}
\label{sv}
 For $ \frp \in \on{NC}_2(2k),\ \# S_n^{\wedge}(\frp) = o(n^{k+1})$.
\end{lemma}
\begin{proof}
Since $\frp$ is noncrossing, we find $\nu \in \ii_{2k-1}$ such
that $\nu \sim_{\frp} \nu + 1$. Suppose that $\nu$ is odd. Then
for any $(\phi, \psi) \in S_n^{\wedge}(\frp)$ we have $\phi(\nu) =
\phi(\nu+1)$. If $\psi(\nu) \neq \psi(\nu+1)$, then by condition
(MP3) there are $o(n^2)$ possibilities for the triplet
$(\phi(\nu), \psi(\nu), \psi(\nu + 1))$, and one sees as in Lemma
\ref{cr} that there are $O(n^{k-1})$ choices for $(\phi, \psi)$
on $\jj := \ii_{2k} \setminus \{\nu, \nu+1\}$. If $\psi(\nu) =
\psi(\nu+1)$, then $\psi(\nu - 1) = \psi(\nu + 2).$ If
$\tilde{\frp}$ is the partition of $\jj$ induced by $\frp$, then
$(\phi|_{\jj}, \psi|_{\jj}) \in S_n^{\wedge}(\tilde{\frp})$. In
this case, hence, $\# S_n^{\wedge}(\frp) \le O(n)\ \#
S_n^{\wedge}(\tilde{\frp}).$ Putting both cases for $\psi(\nu),
\psi(\nu+1)$ together, we obtain
$$ \# S_n^{\wedge}(\frp) \le O(n)\ \# S_n^{\wedge}(\tilde{\frp}) +
o(n^{k+1}).$$ This was proven for $\nu$ odd. In view of the
symmetry of (MP3), the proof for $\nu$ even is analogous, with the roles of
$\phi$ and $\psi$ interchanged. Since
$(\phi, \psi) \in S_n^{\wedge}(\frp)$, iteration of this process
will finally lead to a nearest neighbour pair $\nu \sim_{\frp} \nu
+ 1$ such that $\phi(\nu) \neq \phi(\nu + 1)$ or $\psi(\nu) \neq
\psi(\nu + 1)$. So we end up with
$$ \# S_n^{\wedge}(\frp) \le O(n^{k-1}) o(n^2) + o(n^{k+1}) =
o(n^{k+1}).$$

\end{proof}

Putting the ingredients together, we have that
$$ \frac{1}{n}\ \erw(\on{Tr}((X^* X)^k)) = \frac{1}{n^{k+1}}\
\sum_{\frp \in \calp(2k)}\ \erw \Sigma_n(\frp)$$ is by Corollary
\ref{pair} and Lemmata \ref{ck}, \ref{cr}, and \ref{nc} asymptotically equivalent to
\begin{equation}
\label{qv3} \frac{1}{n^{k+1}}\ \sum_{j=1}^k \sum_{\substack{\frp
\in \on{NC}_2(2k),\\
|\frp \vee \frm| = j}} \erw \Sigma_n(\frp).
\end{equation}
By Lemmata \ref{sv} and \ref{ck} we may replace $\erw
\Sigma_n(\frp)$ by $\erw \Sigma_n^{\vee}(\frp)$ in \eqref{qv3}.
Recall that for $(\phi, \psi) \in S_n^{\vee}(\frp)$, $\phi$ and
$\psi$ are {\it a fortiori} constant on the blocks of $\frp$.
Comparing Lemma \ref{evenodd} with \eqref{qv4}, one sees that this
implies that given $(\phi, \psi) \in  S_n^{\vee}(\frp)$, a block
of $\frp$ corresponds to a matrix element and its complex
conjugate. Invoking \eqref{ssbvar} and Lemma \ref{nc}, we see that
\eqref{qv3} is asymptotically equivalent to
$$ \frac{1}{n^{k+1}}\ \sum_{j=1}^k \sum_{\substack{\frp
\in \on{NC}_2(2k),\\
|\frp \vee \frm| = j}} \sum_{(\phi, \psi) \in S_n^{\vee}(\frp)}\
1\quad  =\frac{1}{n^{k+1}}\ \sum_{j=1}^k \sum_{\substack{\frp
\in \on{NC}_2(2k),\\
|\frp \vee \frm| = j}} s(n)^j\ t(n)^{k - j + 1},
$$
which tends, as $n \to \infty$, to
$$ \sum_{j=1}^k \# \{ \frp
\in \on{NC}_2(2k): |\frp \vee \frm| = j\}\ \kappa^j\ \mu^{k - j +
1} = \sum_{j=1}^k \# \on{NC}^{(j)}(k)\ \kappa^j\ \mu^{k - j +
1},$$ where the last equality follows from Lemma \ref{count}.

\section{The chiral classes}
\label{chiral}

In this section we apply our result about sample covariance matrices to random elements
of the spaces $\mm_n^{\calc}$ from Section \ref{notat}, where $\calc = $ BDI, AIII or CII. It is convenient to consider
the AIII case first. It consists of matrices of the form
$$ \calx_n = \left(
\begin{array}{ll} 0 & X_n \\ X_n^* & 0 \end{array}
\right) \in \cc^{n \times n},$$ with $X_n \in \cc^{s(n) \times
t(n)}$, hence $s(n) + t(n) = n$. We assume that $\lim_{n \to
\infty} \frac{s(n)}{n} = \kappa \in ] 0, 1[$, hence $\lim_{n \to
\infty} \frac{t(n)}{n} = 1 - \kappa =: \mu$. Note that this
framework is more restrictive than the one considered in Section
\ref{sample}. But these restrictions naturally arise if one
considers $X_n$ as a subblock of an element of $\mm_n^{{\rm
AIII}}$, whence $n$ is the natural parameter for asymptotics.
Observe that \begin{equation} \label{qv21}
 \on{Tr}\left(
\begin{array}{ll} 0 & X_n \\ X_n^* & 0 \end{array}
\right)^k = \left\{ \begin{array}{ll} 0 & \text{if $k$ odd}\\
2 \on{Tr} (( X_n^* X_n)^{l}) & \text{if $k = 2l$ even.}\end{array}
\right.
\end{equation}
Assuming that $X_n$ satisfies conditions (MP1), (MP2), (MP3) of Section \ref{sample}, Theorem \ref{thmwish} implies that as $n \to \infty$,
the $2l$-th moment of $\erw L_n(\calx_n)$ converges to
\begin{equation}
\label{momchir}
2 \sum_{j=1}^l \# \on{NC}^{(j)}(l)\ \kappa^j\ (1-\kappa)^{l - j + 1}.
\end{equation}

The special case where the entries of $X_n$ take purely imaginary values yields the same limit for class BDI. Since the extra symmetries of the
CII case are compatible with (MP1), (MP2) and (MP3), we obtain the same limit for this class, too.\\

In \cite{azldp}, in the special case of chiral random matrices with independent Gaussian entries,
the empirical limit distribution $\mu_{\on{ch}, 2}$ of the squared eigenvalues was determined as
\begin{equation}
\label{chir2}
|1 - 2\kappa| \delta_0 + 1_{[a, b]}(x)\  \frac{1}{\pi x} \sqrt{(x -
a) (b - x)}\ dx,
\end{equation}
 where
\begin{equation}
\label{ab}
 a = 1 - 2\sqrt{ \kappa (1 -
\kappa)},\ b = 1 + 2\sqrt{ \kappa (1 - \kappa)}.\end{equation}
Note that \eqref{chir2} differs from the corresponding formula in \cite{azldp},
since a different definition of $L_n$ is used in that paper, and since \eqref{ssbvar} above imposes a condition on complex rather than real variances.
Now, the elegant
approach of Haagerup and Thorbj{\o}rnsen (\cite{HT}) to the moments of the Mar\v{c}enko-Pastur
distribution can be easily adapted to  $\mu_{\on{ch}, 2}$, to the effect that
the $l$-th moment of  $\mu_{\on{ch}, 2}$ is indeed given by \eqref{momchir}. In fact, for $l \ge 1$,
\begin{eqnarray}
\int x^l \mu_{\on{ch},2}(dx) &=& \frac{1}{\pi} \int_{1 - 2\sqrt{\kappa (1 -
\kappa)}}^{1 + 2\sqrt{\kappa (1 - \kappa)}} x^{l - 1} \sqrt{ 4
\kappa (1 - \kappa) - (x - 1)^2}\ dx \nonumber \\
&=& \frac{4 \kappa (1 - \kappa)}{\pi} \int_{- \pi}^0 \sin^2
\theta\ (2 \sqrt{ \kappa ( 1-\kappa)} \cos \theta + 1)^{l - 1}\
d\theta \nonumber \\
&=& \frac{2 \kappa (1 - \kappa)}{\pi} \int_{- \pi}^{\pi} \sin^2
\theta\ (2 \sqrt{ \kappa ( 1-\kappa)} \cos \theta + 1)^{l - 1}\
d\theta. \label{qv20}
\end{eqnarray}
Setting $g(\theta) = (\sqrt{\kappa} + \sqrt{1-\kappa}\ e^{i
\theta})^{l-1}$ and observing that $\sin^2 \theta = \frac12 ( 1 -
\cos 2\theta)$, we see that \eqref{qv20} can be written as
$$ \frac{\kappa (1 - \kappa)}{\pi} \int_{- \pi}^{\pi} \on{Re}(1 - e^{i 2 \theta})\
|g(\theta)|^2 d\theta$$ or as
$$ 2 \kappa (1 - \kappa) \left( \frac{1}{2\pi} \int_{- \pi}^{\pi} |g(\theta)|^2  d\theta\ -
\on{Re} ( \frac{1}{2\pi} \int_{- \pi}^{\pi} h(\theta) \overline{k(\theta)}
d\theta)\right)$$ with $h(\theta) = e^{i\theta} g(\theta),\
k(\theta) = e^{- i\theta} g(\theta).$ Invoking Parseval's formula
and elementary computations with binomial coefficients, we obtain
\begin{eqnarray*}
 \int x^l \mu_{\on{ch},2}(dx) &=& 2 \sum_{j=0}^{l-1} \left\{  {l-1 \choose j}^2
- { l-1 \choose j-1} {l-1 \choose j+1}\right\} \kappa^{j + 1} (1 -
\kappa)^{l - j}\\
 &=& \frac{2}{l} \sum_{j=1}^{l}  {l \choose j} {l \choose j-1} \kappa^{j} (1 -
\kappa)^{l + 1 - j}\\ &=& 2 \sum_{j=1}^{l}  \# (\on{NC}^{(j)}(l))\
\kappa^{j} (1 - \kappa)^{l + 1 - j}.
\end{eqnarray*}
A reference for the last equality is \cite[Cor.\ 9.13]{NSp}. In
view of \eqref{qv21}, we have
\begin{theorem}
If $\calc$ is a chiral class, then
$\erw(L_{\delta n}(X_{\delta n}))$
converges to a probability measure $\mu_{\on{ch}}$ on $\rr$ given by
$$ |1 - 2\kappa| \delta_0 + 1_{[a, b]}(x^2)\  \frac{2}{\pi x} \sqrt{(x^2 -
a) (b - x^2)}\ dx$$ with $a, b$ as in \eqref{ab} above.
\end{theorem}

\bibliographystyle{amsplain}
\bibliography{wignerbiblio}

\end{document}